\documentclass{amsart}
\usepackage[latin1]{inputenc}
\usepackage{amssymb,amsmath}

    {\par\noindent{\scshape Proof:\ }#1}%
    {\mbox{}\hfill$\square$\par\bigskip\par}
\newtheorem{thm}{Theorem}
\newtheorem{defn}[thm]{Definition}
\newtheorem{prop}[thm]{Proposition}

\newtheorem{lem}[thm]{Lemma}

\newcommand{\To}{\longrightarrow}

\begin{document}
\title{Extensions of Boolean isometries}
\author{Antonio Avilés}
\address{Departamento de Matemáticas. Universidad de Murcia. 30100 Murcia, Spain. avileslo@um.es}
\thanks{Author supported by FPU grant of SEEU-MEC, Spain.}

\begin{abstract}
We study when a map between two subsets of a Boolean domain $W$
can be extended to an automorphism of $W$. Under many hypotheses,
if the underlying Boolean algebra is complete or if the sets are
finite or Boolean domains, the necessary and sufficient condition
is that it preserves the Boolean distance between every couple of
points.
\end{abstract}

\maketitle

\section{Introduction}

Boolean domains and Boolean transformations are the Boolean analogues of
algebraic varieties and morphisms of algebraic varieties. We fix once and
for all a Boolean algebra $B$. A Boolean function $f:B^n\To B$ is a
function which admits a polynomial expression in terms of the operations
and elements of $B$, such as for instance $f(x_1,x_2) = (x_1\vee
x_2)\bigtriangleup a$, where $a$ is a fixed element of $B$. A Boolean
domain (over $B$) is a subset $V\subset B^n$ which is the set of solutions
to a Boolean equation, namely
$$ V = \{ (x_1,\ldots,x_n)\in B^n : f(x_1,\ldots,x_n)=0\},$$ for
some Boolean function $f:B^n\To B$. If $U\subset B^n$ and
$V\subset B^m$ are Boolean domains, a map $F:U\To V$ is a Boolean
transformation if there are Boolean functions $F_1,\ldots,
F_m:B^n\To B$ such that $$F(x) = (F_1(x),\ldots,F_m(x))$$ for all
$x\in U$. A Boolean isomorphism is a bijective Boolean
transformation (its inverse map is, in fact, a Boolean
transformation too). Two Boolean domains are isomorphic if there
exists a Boolean isomorphism between them. We must mention the
books \cite{Rud} and~\cite{Ru2} as reference treaties about
Boolean functions and equations.

In this paper, we consider the problem of when a given bijection
between two subsets of a Boolean domain $W$ can be extended to a
Boolean isomorphism from the whole $W$ onto itself. One main
result is the following:

\begin{thm}\label{extension of isomorphisms}
Let $U,V,W\subset B^n$ be Boolean domains with $U\cup V\subset W$
and let $F:U\To V$ be a Boolean isomorphism. Then, $F$ is the
restriction of some Boolean isomorphism $F':W\To W$.
\end{thm}

A Boolean domain $U\subset B^n$ can always be considered as a Boolean
metric space with the metric $d(x,y) = \bigvee_{i=1}^n (x_i\bigtriangleup
y_i)$. A \emph{Boolean metric space} (over $B$) is a set $X$ together with
a symmetric map $d:X\times X\To B$ satisfying the following two
properties: $d(x,y)=0$ if and only if $x=y$, and $d(x,z)\leq d(x,y)\vee
d(y,z)$ for all $x,y,z\in X$. This constitutes a category with maps
$f:X\To Y$ which are \emph{contractive}, that is, $d(f(x),f(y))\leq
d(x,y)$ for all $x,y\in X$. When this inequality is an equality and $f$ is
bijective, then $f$ is called an \emph{isometry}. This concept was early
studied in a series of works like \cite{Blu1}, \cite{Blu2}, \cite{Ell1},
\cite{Ell2} and \cite{Mel1}. In \cite{Avi} the close relation between the
metric and the algebraic structure of Boolean domains, in a more general
context, is investigated. The Boolean transformations between Boolean
domains coincide with the contractive maps and the Boolean isomorphisms
with the isometries. Also, the category of Boolean domains and
transformations is equivalent to the category of CFG-spaces (a subclass of
Boolean metric spaces, whose definition is recalled below) and contractive
maps and therefore Theorem~\ref{extension of isomorphisms} is equivalent
to the following:

\begin{thm}\label{extension of isometries}
Let $U,V,W$ be CFG-spaces with $U\cup V\subset W$ and let $F:U\To
V$ be an isometry. Then, $F$ is the restriction of some isometry
$F':W\To W$.
\end{thm}

A direct consequence of this theorem, together with~\cite[Theorem
1.15]{Avi} is that the necessary and sufficient condition for a
bijection between finite subsets  of a Boolean domain $W$ to be
extended to a Boolean isomorphism of $W$ is to be an isometry
between these two finite sets.

It turns out in fact, that when $B$ is a complete Boolean algebra,
then $U$ and $V$ need not be assumed CFG-spaces:

\begin{thm}\label{complete extension of isometries}
Suppose that $B$ is complete. Let $W$ be a CFG-space, $U,V$ subsets of $W$
and $F:U\To V$ an isometry. Then, $F$ is the restriction of some isometry
$F':W\To W$.
\end{thm}

If $A$ is a $p$-ring for some prime number $p$ (that is, a ring in which
$x^p=x$ and $px=0$ for all $x$) then $A$ happens to be a Boolean metric
space over its ring of idempotents with distance $d(x,y)=(x-y)^{p-1}$.
These spaces were investigated in the papers \cite{Zem} and \cite{Mel2}
which study, among others, problems of extension of isometries. Namely,
\cite[theorem~5]{Zem} is the same statement as our Theorem~\ref{complete
extension of isometries} but only for the particular case in which $W$ is
a $p$-ring.\\

The statement of Theorem~\ref{complete extension of isometries} also holds
for contractive maps instead of isometries:

\begin{thm}\label{complete extension of contractions}
Suppose that $B$ is complete. Let $W$ be a CFG-space, $U,V$
subsets of $W$ and $F:U\To V$ a contractive map. Then, $F$ is the
restriction of some contractive map $F':W\To W$.
\end{thm}

We give examples that the hypothesis of completeness cannot be
weakened in Theorems~\ref{complete extension of isometries} and
\ref{complete extension of contractions}.

\section{Notations}

The operations in Boolean algebras will be denoted as $a \vee b$
and $a\wedge b$ for the supremum and infimum and $a\setminus b$
for the difference, 0 and 1 denote the lowest and greatest
element, $\overline{a}=1\setminus a$ is the complement and
$a\bigtriangleup b = (a\setminus b)\vee (b\setminus a)$ is the
symmetric difference which allows to consider $B$ as a ring with
sum $\bigtriangleup$ and product $\wedge$. Elements
$a_0,\ldots,a_n$ of $B$ are \emph{disjoint} if $a_i\wedge a_j=0$
whenever $i\neq j$ and they are a \emph{partition} if moreover
$a_0\vee\cdots\vee a_n = 1$. The lattice order of $B$ is denoted
as $a\leq b$.

With respect to Boolean metric spaces, the distance will be always
denoted by $d$. The product space of the Boolean metric spaces $X$
and $Y$ is $X\times Y$ with the metric
$$d((x,y),(x',y')) = d(x,x')\vee d(y,y').$$ We will
work in pointed Boolean metric spaces, that is, metric spaces $X$
in which a point $0\in X$ has been fixed. Formally,

\begin{defn}
A pointed Boolean metric space is a couple $(X,0)$ where $X$ is a
Boolean metric space with metric $d$ and $0$ is an element of $X$.
A contractive map between two pointed spaces $f:(X,0)\To (X',0')$
is a contractive map $f:X\To X'$ such that $f(0)=0'$.
\end{defn}

In such spaces we will also use the notation $|x| = d(x,0)$. There is no
deep difference in dealing with pointed spaces but it will be convenient
for technical reasons. We shall make use of several tools in this context,
as \emph{convexity} and \emph{orthogonality}, developed in \cite{Avi},
that are explained below.

Let $a_0,\ldots,a_n$ be a partition of $B$ and $x_0,\ldots,x_n$ be
elements of the metric space $X$. An element $x\in X$ is said to
be a \emph{convex combination} of $x_0,\ldots,x_n$ with
coefficients $a_0,\ldots,a_n$ if $a_i\wedge d(x,x_i)=0$ for all
$i$. In this case we write $a_0x_0+\cdots+a_nx_n=x$.

It turns out that $X$ can be always embedded into a module over
$B$ considered as a ring (sending the fixed element 0 to the zero
of the module) in such a way that these convex combinations
correspond exactly with the usual linear combinations, cf.
\cite[Theorem 1.6]{Avi} and \cite[Proposition 1.11]{Avi}. This
means that the notation is coherent and all the usual properties
of sum and multiplication by scalars apply. When $(X,0)$ is a
pointed metric space then we may suppress the term corresponding
to $0$ in notation $a_00+a_1x_1+\cdots+a_nx_n =
a_1x_1+\cdots+a_nx_n$, where $a_1,\ldots,a_n$ are just disjoint.
We also recall that, in product spaces, convex combinations can be
calculated coordinatewise.\\

A set $S\subset X$ is a \emph{system of generators} of $X$, shortly
$X=conv(S)$, if any element of $X$ can be expressed as a convex
combination of elements of $S$ with some coefficients. We mention the fact
that if two contractive maps coincide on a system of generators, then they
are equal.

A metric space $X$ is a \emph{CFG-space} if it verifies the following two
properties:
\begin{enumerate}
\item It is \emph{convex}, that is, for any $x_0,\ldots,x_n\in X$
and any partition $a_0,\ldots,a_n$ of $B$, the convex combination
$x=a_0x_0+\cdots+a_nx_n$ is an element of $X$. \item It is
\emph{finitely generated}, that is, there is a finite system of
generators of $X$.
\end{enumerate}

We also mention the fact that $X$ is a CFG-space if and only if it
is isometric to a Boolean domain, as it follows from~\cite[Theorem
3.8]{Avi}.

The elements $x$ and $y$ of the pointed space $(X,0)$ are
\emph{orthogonal} ($x\perp y$) if $d(x,y) = |x|\vee |y|$. For a subset
$U\subset (X,0)$ with $0\in U$ we set $$U^\perp = \{y\in X : x\perp y\
\forall x\in U\}.$$

It turns out that $U^\perp$ is a CGF-space provided $U$ is
\cite[Proposition 2.11]{Avi}. The relation of this concept of
orthogonality with the extension of isometries is the following
statement:

\begin{prop}\label{suma de isometrias}
Let $U,X,Y$ be CFG-spaces with $0\in U\subset X$ and $f:(U,0)\To
(Y,0')$ and $g:(U^\perp,0)\To (Y,0')$ be isometries. Then, there
is a unique isometry $f\perp g:(X,0)\To (Y,0')$ which extends both
$f$ and $g$.
\end{prop}

This is the content of Proposition 2.12 in~\cite{Avi} except that there it
is written contractive map instead of isometry. However, it is
straightforward to check in that proof, that if $f$ and $g$ are assumed to
be isometries, then $f\perp g$ that is obtained is again an isometry.

\section{The first extension theorem}

In this section we will prove Theorem~\ref{extension of
isometries}. What we will really prove instead of it will be the
following statement about orthogonal spaces:

\begin{thm}\label{Witt} Let $(X,0)$ be a pointed CFG-space and $U_1,U_2$
CFG-subspaces of $X$ with $0\in U_1\cap U_2$. If $U_1$ is
isometric to $U_2$, then $U_1^\perp$ is isometric to $U_2^\perp$.
\end{thm}

Let us see, first, that Theorem~\ref{extension of isometries} follows from
Theorem~\ref{Witt}. For this, apart from Proposition~\ref{suma de
isometrias}, we need another result \cite[Theorem 4.6]{Avi}, that
CFG-spaces are homogeneous, that is, if $X$ is a CFG-space and $x,y\in X$,
there is an isometry $\phi:X\To X$ such that $\phi(x)=y$. Let $U$, $V$,
$W$ and $F$ be as in the hypotheses of Theorem~\ref{extension of
isometries} and, by homogeneity, fix $0\in U$ and an isometry $\phi:W\To
W$ such that $\phi(F(0))=0$. We apply Theorem~\ref{Witt} to $X=W$,
$U_1=U$, $U_2=\phi(V)$ and we obtain that $U^\perp$ and $\phi(V)^\perp$
are isometric. Again, by homogeneity, we find an isometry
$g:(U^\perp,0)\To(\phi(V)^\perp,0)$. Finally, the map $F' =
\phi^{-1}\circ((\phi\circ F)\perp g)$ is the desired isometry.

Before passing to the proof of Theorem~\ref{Witt}, we must recall
the criteria of isometry and the concept of base developed in
\cite{Avi}.

For a space $X$ and an integer $k>0$, we define an element
$$\alpha_k(X) = \sup\{\bigwedge_{0\leq i<j\leq k}d(u_i,u_j) :
u_0,\ldots,u_k\in X\}$$
 This supremum exists and is indeed attained whenever $X$ is
 either finite or a CFG-space. In the latter case in addition, there
 exists $k_0$ with $\alpha_k(X)=0$ for all $k>k_0$ and
 $\alpha_k(X)\geq \alpha_{k+1}(X)$ for all $k$.
 Another property is that if
 $A$ is a system of generators of $X$, $X=conv(A)$, then
 $\alpha_k(A) = \alpha_k(X)$ for all $k$. The importance of these
 functions is that they determine the isometry classes of
 CFG-spaces: two CFG-spaces $X$ and $Y$ are isometric if and only
 if $\alpha_k(X)=\alpha_k(Y)$ for all $k$, cf.~\cite[\S 4]{Avi}.

 Another result that we need is the existence of \emph{bases}:
 Any pointed CFG-space $(X,0)$ has a base, that is, a set
 $\{x_1,\ldots,x_n\}$ such that\begin{enumerate}
\item $X = conv(0,x_1,\ldots,x_n)$,
  \item $x_i\perp x_j$ for any
 $i\neq j$,
 \item $\alpha_i(X)=|x_i|>0$ for $i=1,\ldots,n$.
\end{enumerate}

We point out that condition (1) above implies that $\alpha_i(X) =
0$ for $i>n$. The following lemma investigates the relation
between the functions $\alpha_k(U)$, $\alpha_k(U^\perp)$ and
$\alpha_k(X)$ when $U$ is a CFG-subspace of $X$. It will be useful
now to convene that $\alpha_0(Y) = 1$ for any space $Y$.

\begin{lem}\label{amplitudes de la suma} Let $(X,0)$ be a CFG-space and $U$ a CFG-subspace with
$0\in U$. Then, for all $n\in\mathbb{N}$,
$$\alpha_{n}(X) = \bigvee_{i=0}^{n}\alpha_{i}(U)\wedge\alpha_{n-i}(U^{\perp}).$$
\end{lem}

PROOF: Take bases $B_1 = \{x_1,\ldots,x_r\}$ and $B_2 =
\{y_1,\ldots,y_s\}$ of $(U,0)$ and $(U^\perp,0)$, respectively and define
$B= B_1\cup B_2\cup\{0\}$. From \cite[Proposition 2.11]{Avi} we have
$X=conv(U\cup U^\perp)$ and hence $ X = conv(B)$ and $\alpha_n(X) =
\alpha_n(B)$. Now the result follows by applying the definition of the
function $\alpha_n$ to that set, having in
mind the relations \begin{eqnarray} |x_1|\geq|x_2|\geq\cdots,\\
|y_1|\geq|y_2|\geq\cdots,\\ d(x_i,x_j)=|x_i|\vee|x_j| = |x_{\min(i,j)}|,\\
d(y_i,y_j)=|y_i|\vee|y_j|=|y_{\min(i,j)}|.
\end{eqnarray} Namely, for a subset $A$ of $B$ we define
$$\phi(A) = \bigwedge_{u,v\in A, u\neq v}d(u,v),$$ so that $\alpha_n(B)$ is
the supremum of all $\phi(A)$ when $A$ runs over all subsets of
$B$ of cardinality $n+1$. Whenever $n-s\leq i\leq r$, we can
consider the set
$$A_i=\{0,x_1,\ldots,x_i,y_1,\ldots,y_{n-i}\}$$
of cardinality $n+1$, so that $\alpha_n(X)\geq \phi(A_i)$ and by the
relations mentioned above, it is easily calculated that $\phi(A_i) =
|x_i|\wedge |y_{n-i}| = \alpha_i(U)\wedge\alpha_{n-i}(U^\perp)$. When
$n-s\leq i\leq r$ does not hold, then
$\alpha_{i}(U)\wedge\alpha_{n-i}(U^{\perp}) = 0$. This proves that
$\alpha_{n}(X) \geq
\bigvee_{i=0}^{n}\alpha_{i}(U)\wedge\alpha_{n-i}(U^{\perp}).$ For the
other inequality, we take an arbitrary subset $A$ of $B$ of cardinality
$n+1$ and we shall prove that $\phi(A)\leq
\bigvee_{i=0}^{n}\alpha_{i}(U)\wedge\alpha_{n-i}(U^{\perp})$. For such an
$A$, we find $i_1<\cdots<i_t$ and $j_1<\cdots<j_u$ such that
\begin{eqnarray*} A\cap B_1 &=& \{x_{i_1},\ldots,x_{i_t}\},\\
A\cap B_2 &=& \{y_{j_1},\ldots,y_{j_u}\}.
\end{eqnarray*}

Now, if $0\in A$ then $t+u = n$ and using relations $(1) - (4)$
above
$$\phi(A) \leq d(0,x_{i_t})\wedge d(0,y_{j_u}) =  |x_{i_t}|\wedge
|y_{j_u}|\leq |x_t|\wedge |y_{u}| = \alpha_t(U)\wedge\alpha_u(U^\perp).$$
On the other hand, if $0\not\in A$, then $u+t=n+1$ and calculating again,
\begin{eqnarray*}\text{if }t,u\geq 2,\ \ \phi(A) &\leq& d(x_{i_t},y_{j_u})\wedge
d(x_{i_t},x_{i_{t-1}})\wedge d(y_{j_u},y_{j_{u-1}})\\
&=& (|x_{i_t}|\vee |y_{j_u}|)\wedge |x_{i_{t-1}}|\wedge
|y_{j_{u-1}}|\\
&=& (|x_{i_t}|\wedge|y_{j_{u-1}}|)\vee (|x_{i_{t-1}}|\wedge
|y_{j_u}|)\\
 &\leq& (|x_{t}|\wedge|y_{u-1}|)\vee (|x_{t-1}|\wedge
|y_{u}|),
\end{eqnarray*}

\begin{eqnarray*}\text{if }t=1, u>1,\ \ \phi(A) &\leq& d(x_{i_1},y_{j_u})\wedge
 d(y_{j_u},y_{j_{u-1}})\\
&=& (|x_{i_1}|\vee |y_{j_u}|)\wedge
|y_{j_{u-1}}|\\
&=& (|x_{i_1}|\wedge|y_{j_{u-1}}|)\vee |y_{j_u}|)\\
 &\leq& (|x_{1}|\wedge|y_{u-1}|)\vee |y_{u}|\\ &=&
 (\alpha_1(U)\wedge \alpha_{u-1}(U^\perp))\vee (\alpha_0(U)\wedge
 \alpha_u(U^\perp)),
\end{eqnarray*}
 and the other cases are checked similarly.$\qed$\\

PROOF OF THEOREM~\ref{Witt}: For every $i\in\mathbb{N}$, we set
$$a_i=\alpha_i(U_1)=\alpha_i(U_2),\ b_i=\alpha_i(X),\
r_i=\alpha_i(U_1^\perp),\ s_i=\alpha_i(U_2^\perp).$$
 What we must prove is that $r_i=s_i$ for every $i$. Let $d$ be
 the greatest integer with $\alpha_{d}(X)>0$. Clearly, $r_i=s_i=0$
 for all $i>d$ and by Lemma~\ref{amplitudes de la suma} both
 $(r_i)_{i=1}^d$ and $(s_i)_{i=1}^d$ are solutions to the
 following system of equations in the variables $x_1,\ldots,x_d$:
 \begin{eqnarray}
 x_1\geq\cdots\geq x_d, & &\\
 \bigvee_{i=0}^n (a_{n-i}\wedge x_i) = b_n & & \ \ \
 n=1,\ldots,d+1,
 \end{eqnarray}
where $x_0=a_0=1$ and $x_{d+1} = 0$ are constants.

Hence, we must see that this system of equations has a unique solution,
under the hypotheses that $b_1\geq\cdots\geq b_{d+1}=0$,
$a_1\geq\cdots\geq a_{d+1}=0$ and $a_i\leq b_i$ for all $i$. We need,
therefore, a criterion to ensure the uniqueness of solutions of a certain
system of Boolean equations, which is provided by the following lemma:

\begin{lem}\label{unicidad de soluciones}
Let $Y$ be a CFG-space and $\{y_0,\ldots,y_n\}$ a system of
generators of $Y$ such that $d(y_i,y_j)=1$ for all $i\neq j$. Let
$f:Y\To B$ be a contractive function such that $f(y_i)\vee
f(y_j)=1$ for all $i\neq j$. If the equation $f(x)=0$ has a
solution for $x\in Y$, then this solution is unique.
\end{lem}

PROOF: Notice that, even if $i=j$ we always have $d(y_i,y_j)\leq
f(y_i)\vee f(y_j)$ for all $i,j=0,\ldots,n$. The set of all
couples $(y_i,y_j)$ is a system of generators of the product space
$Y\times Y$. We consider the function $h(x,y) = d(x,y)\setminus
(f(x)\vee f(y))$ on $Y\times Y$. First, we notice that $h$ is
contractive. The map $(x,y)\mapsto f(x)\vee f(y)$ is contractive
since it is the composition of contractive maps $(x,y)\mapsto
(f(x),f(y))$ and $(a,b)\mapsto a\vee b$. The map $(x,y)\mapsto
d(x,y)$ is also contractive, cf. property (3') after
\cite[Definition 1.1]{Avi}. Hence $h$ is contractive since it is a
Boolean operation of two contractive maps. On the other hand, $h$
is equal to zero on the system of generators $\{(y_i,y_j)\}$ and
therefore, it is constant equal to zero on all $Y\times Y$. Hence,
if $f(x)=f(y)=0$, then $d(x,y)=0$ and $x=y$.$\qed$\\

Back to the proof of Theorem~\ref{Witt}, we shall apply
Lemma~\ref{unicidad de soluciones} to $$Y=\{(x_1,\ldots,x_d)\in
B^d : x_1\geq\cdots\geq x_d\},$$
 which is a metric space with the usual metric
 $d(x,x')=\bigvee_{i=1}^{d}(x_i\bigtriangleup x'_i)$. It is
 checked in \cite{Avi} that in these metric spaces, convex
 combinations are calculated simply coordinatewise in the natural
 way. It is straightforward to check that in fact, $Y$ is a
 CFG-space with the set of generators

$\begin{array}{r}
   y_0 = (0,0,\ldots,0,0),\\
  y_1 = (1,0,\ldots,0,0), \\
   \cdots\\
   y_{d-1} = (1,1,\ldots,1,0),\\
   y_d = (1,1,\ldots,1,1).\\
\end{array}$

Namely, if $c = (c_1,\cdots,c_d)$ then $c= (c_1\setminus c_2) y_1 +
(c_2\setminus c_3)y_2 + \cdots + c_d y_d + (1\setminus\bigvee c_i)y_0$.
After~\cite[Theorem 3.8]{Avi}, the contractive functions from $Y$ to $B$
are exactly the Boolean functions. We will finish the proof provided we
can apply Lemma~\ref{unicidad de soluciones} to the Boolean function
$f(x)=\bigvee_{n=1}^{d+1}f_n(x)$, where

$$f_n(x_1,\ldots,x_d) = b_n\bigtriangleup
\left(\bigvee_{i=0}^n a_{n-i}\wedge x_i\right).$$

It remains to check that $f(y_j)\vee f(y_k) = 1$ whenever
$j,k=0,\ldots,d$, $j\neq k$. First, we calculate the value of the
$f_n(y_j)$'s. For notational simplicity we convene that
$(y_j)_0=1$.\begin{eqnarray*}
    f_{n}(y_{j}) &=& b_{n}\bigtriangleup \bigvee_{i=0}^{n}a_{n-i}\wedge(y_{j})_{i}
    = b_{n}\bigtriangleup (a_{n-0}\vee a_{n-1}\vee\cdots\vee a_{n-j})\\ &=&
    a_{n-j}\bigtriangleup b_{n} \text{ if }j<n;\\
    & & \\
    f_{n}(y_{j}) &=& b_{n}\bigtriangleup \bigvee_{i=0}^{n}a_{n-i}\wedge(y_{j})_{i}
    = b_{n}\bigtriangleup a_{0}= b_{n}\bigtriangleup 1\\ &=& \overline{b_{n}} \text{ if }j\geq
    n.
    \end{eqnarray*}

The value of the $f(y_j)$'s is then

\begin{eqnarray*}
        f(y_{0}) &=& (a_{1}\bigtriangleup b_{1})\vee (a_{2}\bigtriangleup b_{2})\vee\cdots\vee
        (a_{d}\bigtriangleup b_{d})\vee 0; \\
        f(y_{1}) &=& \overline{b_{1}}\vee (a_{1}\bigtriangleup b_{2})\vee\cdots\vee
        (a_{d-1}\bigtriangleup b_{d})\vee a_{d}; \\
        f(y_{2}) &=& \overline{b_{1}}\vee \overline{b_{2}}\vee
        (a_{1}\bigtriangleup b_{3})\vee\cdots\vee (a_{d-2}\bigtriangleup b_{d})\vee a_{d-1}; \\
        & & \cdots \\
        f(y_{j}) &=&
        \overline{b_{1}}\vee\cdots\overline{b_{j}}\vee
        (a_{1}\bigtriangleup b_{j+1})\vee\cdots\vee(a_{d-j}\bigtriangleup b_{d})\vee
        a_{d-j+ 1};\\
        & & \cdots \\
        f(y_{d-1}) &=& \overline{b_{1}}\vee\cdots\vee \overline{b_{d-1}}\vee
        (a_{1}\bigtriangleup b_{d})\vee a_{2}; \\
        f(y_{d}) &=& \overline{b_{1}}\vee\cdots\vee \overline{b_{d}}\vee
        a_{1}.
      \end{eqnarray*}

We can simplify since $b_1\geq b_2\geq\cdots\geq b_d$:

\begin{eqnarray*}
        f(y_{0}) &=& (a_{1}\bigtriangleup b_{1})\vee (a_{2}\bigtriangleup b_{2})\vee\cdots\vee
        (a_{d}\bigtriangleup b_{d}); \\
        f(y_{1}) &=& \overline{b_{1}}\vee (a_{1}\bigtriangleup b_{2})\vee\cdots\vee
        (a_{d-1}\bigtriangleup b_{d})\vee a_{d}; \\
        f(y_{2}) &=& \overline{b_{2}}\vee
        (a_{1}\bigtriangleup b_{3})\vee\cdots\vee (a_{d-2}\bigtriangleup b_{d})\vee a_{d-1}; \\
        & & \cdots \\
        f(y_{j}) &=&
       \overline{b_{j}}\vee
        (a_{1}\bigtriangleup b_{j+1})\vee\cdots\vee(a_{d-j}\bigtriangleup b_{d})\vee
        a_{d-j+ 1};\\
        & & \cdots \\
        f(y_{d-1}) &=& \overline{b_{d-1}}\vee (a_{1}\bigtriangleup b_{d})\vee a_{2}; \\
        f(y_{d}) &=& \overline{b_{d}}\vee a_{1}.
      \end{eqnarray*}

Now, we fix $i$, $j$ and $a_1\geq\cdots\geq a_d$. We must see that
for any $(b_1\ldots,b_d)\in Y$, $f(y_i)\vee f(y_j)=1$. Again, the
function $\phi(b)$ which associates to each $b=(b_1,\ldots,b_d)\in
Y$ the corresponding value of $\phi(b) = f(y_i)\vee f(y_j)$ is a
Boolean function, and in order to see that $\phi$ is constant
equal to one on $Y$ it is enough to check that $\phi(y_k)=1$ for
$k=0,\ldots,d$. For $(b_1,\ldots,b_d) = y_k$ we obtain:

\begin{eqnarray*}
        f(y_{0}) &=& (a_{1}\bigtriangleup 1)\vee\dots\vee (a_{k}\bigtriangleup 1)\vee a_{k+1}\vee\cdots
        \vee a_{d}; \\
        f(y_{1}) &=& (a_{1}\bigtriangleup 1)\vee\cdots\vee (a_{k-1}\bigtriangleup 1)\vee a_{k} \vee\cdots
        \vee a_{d}; \\
        & & \cdots \\
        f(y_{j}) &=& (a_{1}\bigtriangleup 1)\vee\cdots\vee (a_{k-j}\bigtriangleup 1)\vee
        a_{k-j+1}\vee\cdots\vee a_{d-j+1};\\
        & & \cdots\\
        f(y_{k-1}) &=& (a_{1}\bigtriangleup 1)\vee a_{2} \vee\cdots\vee a_{d-k+2}; \\
        f(y_{k}) &=& a_{1}\vee a_{2} \vee\cdots\vee a_{d-k+1}; \\
        f(y_{k+1}) &=& f(y_{k+2}) = \cdots = f(y_{d}) = 1.
      \end{eqnarray*}

Now, it is clear that $f(y_i)\vee f(y_j)=1$ for $i\neq j$ because
if $i<j$ then $a_{k-j+1}\leq f(y_j)$ and $a_{k-j+1}\bigtriangleup
1\leq f(y_i)$. This finishes the proof of Theorem~\ref{Witt} and
hence, also the proofs of Theorems~\ref{extension of isometries}
and~\ref{extension of isomorphisms}.$\qed$\\

\section{The second extension theorem}

In this section we prove Theorems~\ref{complete extension of
isometries} and~\ref{complete extension of contractions}. Hence,
we assume from now on that our fixed Boolean algebra $B$ is
complete, that is, that whenever $S$ is a subset of $B$ there
exists $s=\bigvee S\in B$ the supremum of $S$. We recall that the
distributivity law still holds in the infinite case:
$x\wedge\bigvee\{y_i : i\in I\} = \bigvee\{x\wedge y_i : i\in I\}$
whenever $x\in B$ and $y_i\in B$ for all $i\in I$.

\begin{lem}\label{supremo de contractivas}
Let $X$ be a metric space over $B$ and $\{f_i:X\To B\}_{i\in I}$ a
family of contractive maps. Then, the pointwise supremum $\bigvee
f_i$ is again a contractive map.
\end{lem}

PROOF: Recall that the metric on $B$ is given by $d(x,y)=x\bigtriangleup
y$ and hence $f:X\To B$ is contractive if and only if $f(x)\bigtriangleup
f(y)\leq d(x,y)$ for all $x,y\in X$. Moreover, this can be rewritten as
$$\overline{d(x,y)}\wedge f(y)\leq f(x)\leq f(y)\vee d(x,y)$$
for all $x,y\in X$. With this characterization and using the
infinite distributivity law, the proof of the lemma becomes
apparent.$\qed$

\begin{lem}\label{interseccion de CFGs}
Let $X$ be a CFG space over the complete Boolean algebra $B$ and
let $\{K_i\}_{i\in I}$ be a family of CFG-subspaces of $X$. Then
$\bigcap_{I}K_i$ is a CFG-space.
\end{lem}

PROOF: By~\cite[Lemma~3.5]{Avi} a subspace $K\subset X$ is a
CFG-space if and only if there exists $f:X\To B$ contractive with
$K=f^{-1}(\{0\})$. This together with Lemma~\ref{supremo de
contractivas} proves the Lemma.$\qed$

By Lemma~\ref{interseccion de CFGs}, given a subset $U$ of a CFG-space
$X$, we can consider $Conv(U)$ the least CFG-space that contains $U$,
obtained as the intersection of all CFG-subspaces that contain $U$. Any
nonprincipal $I$ ideal of $B$ is an example in which $I=conv(I)\neq
Conv(I)$ since $I$ is convex but not a CFG-space.

\begin{thm}\label{extension a Conv}
Let $X$ and $Y$ be CFG-spaces over the complete Boolean algebra
$B$ and let $f:U\To V$ be a contractive map between two arbitrary
subsets $U\subset X$ and $V\subset Y$. Then there is a unique
contractive map $Conv(f):Conv(U)\To Conv(V)$ that extends $f$. In
addition, if $f$ is an isometry, so is $Conv(f)$.
\end{thm}

Notice that Theorem~\ref{complete extension of isometries} is a
direct consequence of Theorem~\ref{extension a Conv} above
together with Theorem~\ref{extension of isometries}, while
Theorem~\ref{complete extension of contractions} follows from
Theorem~\ref{extension a Conv} and \cite[Proposition~2.12]{Avi}.\\

PROOF OF THEOREM~\ref{extension a Conv}: First, we check that
$Conv(f)$, provided it exists, is uniquely determined. Suppose
that $g,h:Conv(U)\To Conv(V)$ are two contractive extensions of
$f$. Then the set
$$ K = \{x\in Conv(U) : d(g(x),h(x))=0\}$$
is, by~\cite[Lemma~3.5]{Avi} a CFG-space which contains $U$, hence
$Conv(U)\subset K$ and $g=h$.

 For the existence of $Conv(f)$, we prove first a
particular case, namely, that any contractive function $f:U\To B$
extends to a contractive map $G:Conv(U)\To B$. For every $u\in U$
we consider the contractive map $g_u:Conv(U)\To B$ given by
$$g_u(x) = f(u)\setminus d(u,x)$$
and we set $G=\bigvee\{g_u : u\in U\}$. On the one hand, for any $u\in U$,
$f(u)= g_u(u)\leq G(u)$. On the other hand for any $u,v\in U$,
$f(u)\bigtriangleup f(v)\leq d(u,v)$ and hence $f(v)\geq f(u)\setminus
d(u,v)=g_u(v)$, so taking suprema over $U$, also $f(v)\geq G(v)$. Now we
pass to the general case and we use the fact that $Y$ can be viewed as a
subspace of $B^n$ for some natural number $n$. Extending coordinate by
coordinate, we know that there is a contractive map $h:Conv(U)\To B^n$
which extends $f$. It remains to show that the range of $h$ verifies
$h(Conv(U))\subset Conv(V)\subset Y$. Again, by~\cite[Lemma~3.5]{Avi}
there is a contractive map $s:B^n\To B$ such that $Conv(V)=s^{-1}(\{0\})$.
Notice that for every $u\in U$, $h(u)\in V\subset Conv(V)=s^{-1}(\{0\})$
so $s(h(u))=0$. Therefore the composed map $s\circ h:Conv(U)\To B$ is a
contractive map which extends the constant map $c:U\To B$, $c(u)=0$. By
the uniqueness of extensions to $Conv(U)$ that we have already proved, we
obtain that $s\circ h=0$, so $h(Conv(U))\subset s^{-1}(\{0\})=Conv(V)$.

With respect to the last assertion of the theorem, if $f$ is an
isometry then $f^{-1}:V\To U$ is a contractive map and
$Conv(f^{-1})$ must be a contractive inverse map for $Conv(f)$
(since the compositions in both senses are contractive extensions
of the identity maps in $Conv(U)$ and $Conv(V)$). This implies
that $Conv(f)$ is an isometry.$\qed$

We finish by presenting an example which shows that the hypotheses
of Theorems~\ref{extension of isometries} and \ref{complete
extension of isometries} cannot be essentially weakened.

Assuming that $B$ is not complete we construct a CFG space $X$ and an
isometry $f:U\To V$ between subsets of $X$ which cannot be extended to any
contractive map $F:X\To X$. Take $S$ a subset of $B$ which does not have a
supremum and set $$ I = \{a\in B : \exists a_1,\ldots,a_n\in S : a\leq
a_1\vee\cdots\vee a_n\};$$ the ideal generated by $S$ which neither has a
supremum. Namely, if $x$ were the supremum of $I$, then it would be also
the supremum of $S$ because $S$ and $I$ have the same upper bounds: if $y$
is an upper bound of $S$ and $a\in I$, then $a\leq a_1\vee\cdots\vee a_n$
for some elements $a_i\in S$, so that $a_i\leq y$ for all $i$ and finally
$a\leq y$. Set
\begin{eqnarray*} J &=& \{a\in B : a\wedge x =0\ \forall x\in I\},\\
I+J &=& \{a\bigtriangleup b : a\in I, b\in J\},\\ X &=& \{(x,y)\in
B^2 : x\wedge y = 0\},\\ V &=& \{(x,y)\in X : x\in I, y\in J\},\\
U &=& \{(z,0)\in X : z\in I+J\}.\end{eqnarray*} Observe that $X$ is a
CFG-space since it is a Boolean domain, in fact $X= conv\{(0,0), (0,1),
(1,0)\}$. The isometry is $f=g^{-1}$, the inverse map of $g:V\To U$ given
by $g(x,y)=(x\bigtriangleup y,0)$. Namely $g$ is an isometry because it is
clearly onto and for any $x,x'\in I$ and $y,y'\in J$,
\begin{eqnarray*} d(g(x,y),g(x',y')) &=&
x\bigtriangleup y\bigtriangleup x'\bigtriangleup y'=(x\bigtriangleup x')\bigtriangleup (y\bigtriangleup y');\\
d((x,y),(x',y')) &=& (x\bigtriangleup x')\vee (y\bigtriangleup
y')\end{eqnarray*} and the two expressions are equal because
$x\bigtriangleup x'\in I$ and $y\bigtriangleup y'\in J$, so they are
disjoint.

Suppose that we could extend $f$ to some contractive map $F:X\To X$. We
claim that if $F(1,0)=(a,b)$ then $a$ is the supremum of $I$, which is a
contradiction. Namely, for every $x\in I$,
$$ (x\bigtriangleup a)\vee b = d((x,0),(a,b)) = d(F(x,0),F(1,0))\leq d((x,0),(1,0))=\overline{x} $$
so that $x\leq a$ and analogously for every $y\in J$,
$$ (y\bigtriangleup b)\vee a = d((0,y),(a,b)) = d(F(y,0),F(1,0))\leq d((y,0),(1,0))=\overline{y} $$
and $y\leq b$. This means that $a$ is an upper bound of $I$ and $b$ an
upper bound of $J$. If $c$ is now an arbitrary upper bound of $I$ then
$\overline{c}\in J$, so $\overline{c}\leq b$, so $a\wedge\overline{c}\leq
a\wedge b=0$ and $a\leq c$.

Observe that the space $X$ in the example is ``two-dimensional''. In fact
the case $X=B$ is special and even if $B$ is not complete, arbitrary
isometries between subsets can be always extended. This is because if
$f:U\To V$ is an isometry between $U,V\subset B$ then $f(x)\bigtriangleup
f(y) = x\bigtriangleup y$ for all $x,y\in U$ and this implies that the
function $x\bigtriangleup f(x)$ is constant equal to some $a\in B$, and
then $F(x)=a\bigtriangleup x$ is an isometry of $B$ that extends $F$.
However, this particularity does not apply when we consider extensions of
contractive maps instead of isometries. Take for instance two infinite
sets $M\subset\Omega$ and $B$ the Boolean algebra of the finite or
cofinite subsets of $\Omega$ and $U\subset B$ the family of the finite
subsets of $\Omega$. Then the contractive map $f:U\To U$ given by
$f(x)=M\cap x$ cannot be contractively extended to $B$.

\end{document}